\documentclass[11pt]{article}
\usepackage{amssymb}
\usepackage{amsmath}
\newtheorem{precor}{{\bf Corollary}}

\newenvironment{cor}{\begin{precor}{\hspace{-0.5
               em}{\bf.\ }}}{\end{precor}}
\newtheorem{precon}{{\bf Conjecture}}

\newtheorem{predefin}{{\bf Definition}}

\newenvironment{defin}[1]{\begin{predefin}{\hspace{-0.5
                   em}{\bf.\ }}{\rm #1}\hfill{$\spadesuit$}}{\end{predefin}}
\newtheorem{preexm}{{\bf Example}}

\newtheorem{preappl}{{\bf Application}}

\newtheorem{prelem}{{\bf Lemma}}

\newenvironment{lem}{\begin{prelem}{\hspace{-0.5
               em}{\bf.\ }}}{\end{prelem}}
\newtheorem{preproof}{{\bf Proof.\ }}

\newenvironment{proof}[1]{\begin{preproof}{\rm
               #1}\hfill{$\blacksquare$}}{\end{preproof}}
\newtheorem{prethm}{{\bf Theorem}}

\newenvironment{thm}{\begin{prethm}{\hspace{-0.5
               em}{\bf.\ }}}{\end{prethm}}
\newtheorem{prealphthm}{{\bf Theorem}}

\newenvironment{alphthm}{\begin{prealphthm}{\hspace{-0.5
               em}{\bf.\ }}}{\end{prealphthm}}
\newtheorem{prealphlem}{{\bf Lemma}}

\newenvironment{alphlem}{\begin{prealphlem}{\hspace{-0.5
               em}{\bf.\ }}}{\end{prealphlem}}
\newtheorem{prepro}{{\bf Proposition}}

\newtheorem{prequ}{{\bf Question}}

\newtheorem{prealphqu}{{\bf Question}}

\newenvironment{alphqu}{\begin{prealphqu}{\hspace{-0.5
               em}{\bf.\ }}}{\end{prealphqu}}
\newtheorem{preprb}{{\bf Problem}}

\def\conct[#1,#2]{\mbox {${#1} \leftrightarrow {#2}$}}
\def\dconct[#1,#2]{\mbox {${#1} \rightarrow {#2}$}}
\def\deg[#1,#2]{\mbox {$d_{_{#1}}(#2)$}}
\def\mindeg[#1]{\mbox {$\delta_{_{#1}}$}}
\def\maxdeg[#1]{\mbox {$\Delta_{_{#1}}$}}
\def\outdeg[#1,#2]{\mbox {$d_{_{#1}}^{^+}(#2)$}}
\def\minoutdeg[#1]{\mbox {$\delta_{_{#1}}^{^+}$}}
\def\maxoutdeg[#1]{\mbox {$\Delta_{_{#1}}^{^+}$}}
\def\indeg[#1,#2]{\mbox {$d_{_{#1}}^{^-}(#2)$}}
\def\minindeg[#1]{\mbox {$\delta_{_{#1}}^{^-}$}}
\def\maxindeg[#1]{\mbox {$\Delta_{_{#1}}^{^-}$}}

\def\dre[#1,#2,#3]{\mbox {${\cal E}_{_{#3}}(#1,#2)$}}
\def\var[#1,#2]{\mbox {${\rm Var}_{_{#1}}(#2)$}}
\def\ls[#1]{\mbox {$\xi^{^{#1}}$}}
\def\hom[#1,#2]{\mbox {${\rm Hom}({#1},{#2})$}}
\def\onvhom[#1,#2]{\mbox {${\rm Hom^{v}}(#1,#2)$}}
\def\onehom[#1,#2]{\mbox {${\rm Hom^{e}}(#1,#2)$}}
\def\core[#1]{\mbox {$#1^{^{\bullet}}$}}
\def\cay[#1,#2]{\mbox {${\rm Cay}({#1},{#2})$}}
\def\cays[#1,#2]{\mbox {${\rm Cay_{s}}({#1},{#2})$}}
\def\dirc[#1]{\mbox {$\stackrel{\rightarrow}{C}_{_{#1}}$}}
\def\cycl[#1]{\mbox {${\bf Z}_{_{#1}}$}}

\date{}
\begin{document}
\footnotetext[1]{$\ast$This research was in part supported by
Shahid Beheshti University.}
\begin{center}
{\Large \bf On Coloring Properties of Graph Powers}\\
\vspace*{0.5cm}
{\bf Hossein Hajiabolhassan$^\ast$}\\
{\it Department of Mathematical Sciences}\\
{\it Shahid Beheshti University,  G.C.,}\\
{\it P.O. Box {\rm 19839-63113}, Tehran, Iran}\\
{hhaji@sbu.ac.ir}\\
{\bf Ali Taherkhani}\\
{\it Department of Mathematics}\\
{\it Institute for Advanced Studies in Basic Sciences }\\
{\it P.O. Box {\rm 45195-1159}, Zanjan {\rm 45195}, Iran}\\
{\tt ali.taherkhani@iasbs.ac.ir}\\
\end{center}
\begin{abstract}
\noindent This paper studies some coloring properties of graph powers.
We  show that $\chi_c(G^{^{\frac{2r+1}{2s+1}}})=\frac{(2s+1)\chi_c(G)}{(s-r)\chi_c(G)+2r+1}$
provided that $\chi_c(G^{^{\frac{2r+1}{2s+1}}})< 4$.  As a consequence,
one can see that if ${2r+1 \over 2s+1} \leq {\chi_c(G) \over 3(\chi_c(G)-2)}$,
then $\chi_c(G^{^{\frac{2r+1}{2s+1}}})=\frac{(2s+1)\chi_c(G)}{(s-r)\chi_c(G)+2r+1}$.
 In particular, $\chi_c(K_{3n+1}^{^{1\over3}})={9n+3\over 3n+2}$ and $K_{3n+1}^{^{1\over3}}$
 has no subgraph with circular chromatic number equal to ${6n+1\over 2n+1}$. This
 provides a negative answer to a question asked  in
[Xuding Zhu, Circular chromatic number: a survey, Discrete Math.,
229(1-3):371--410, 2001]. Also, we present
 an upper bound for the fractional chromatic number of subdivision graphs. Precisely, we show that
 $\chi_f(G^{^{\frac{1}{2s+1}}})\leq \frac{(2s+1)\chi_f(G)}{s\chi_f(G)+1}$. Finally,
 we investigate the $n$th multichromatic number of subdivision graphs.
\\

\end{abstract}
\section{Introduction}
It was shown in \cite{larsen} that one can compute the fractional
chromatic number of $M(G)$ in terms of that of $G$, where $M(G)$
stands for the Mycielskian of $G$. There are a few interesting
and similar results for the circular chromatic number. Hence, it
is of interest to find a map or a functor ${\cal F}$ from the
category of graphs to itself such that, for any graph $G$, it is
possible to determine the exact value of the circular chromatic
number of ${\cal F}(G)$ in terms of that of $G$. In this paper, we
show that graph powers can be considered as such fanctors (graph
powers preserve the graph homomorphism).

In Section~$1$, we set up notation and terminology. Section~$2$
establishes the tight relation between the circular chromatic
number and graph powers. In fact, we show that it is possible to
determine the circular chromatic number of $G^{2r+1\over 2s+1}$
in terms of that of $G$ provided that ${2r+1\over 2s+1}$ is
sufficiently small. In Section~$3$, we investigate the fractional
chromatic number and the $n$th multichromatic number of subdivison
graphs.

Throughout this paper we consider finite simple graphs which have
no loops and multiple edges. For a given graph $G$, the notation
${\rm og}(G)$ stands for the odd girth of $G$. We denote by $[m]$
the set $\{1,2,\ldots,m\}$. Let $G$ and $H$ be two graphs. A
homomorphism from $G$ to $H$ is a mapping $f:V(G)\longrightarrow
V(H)$ such that $f(u)f(v)\in E(H)$ whenever $uv\in V(G)$. We
write $G\longrightarrow H$ if there exists a homomorphism from
$G$ to $H$. Two graphs $G$ and $H$ are homomorphically equivalent
if $G\longrightarrow H$ and $H\longrightarrow G$ and it is
indicated by the symbol $G\longleftrightarrow H$.


Let $d$ and $n$ be positive integers, where $n\geq 2d$. The
circular complete graph $K_{n\over d}$ has the vertex set
$\{0,1,\ldots,n-1\}$ in which $ij$ is an edge if and only if
$d\leq |i-j|\leq n-d$.
 An $(n,d)-$coloring of graph $G$ is a homomorphism from $G$ to the circular complete graph $K_{n\over d}$.
 The circular chromatic number $\chi_c(G)$ of $G$ is defined as
 $$\chi_c(G)=\inf\{{n\over d}| G\,\, {\rm admits\,\, an}\,\, (n,d)-{\rm coloring}\}.$$

Two kinds of graph powers were introduced  in
\cite{hajiabolhassan, MR2587748}. Especially, it was illustrated
that there is a tight relationship between graph powers and the
circular chromatic number. Also, the connection between graph
homomorphism and graph powers has been studied  in
\cite{hajiabolhassan, MR2587748, MR2171371}.

For a graph $G$, let $G^{^{k}}$ be the $k$th power of $G$, which
is obtained on the vertex set $V(G)$, by connecting any two
vertices $u$ and $v$ for which there exists a walk of length $k$
between $u$ and $v$ in $G$. Also, assume that $G^{^{1\over s}}$
is the graph obtained by replacing each edge of $G$ with the path
$P_{s+1}$. Set $G^{^{r\over s}}= (G^{^{1\over s}})^r$. This
power, called fractional power as a functor, preserves the graph
homomorphism.
 In this terminology, we have the following lemma.
 \begin{alphlem}\label{hada}{\rm\cite{MR2587748}}
Let $r$ and $s$ be positive integers and G be a  graph. Then
$$  \begin{array}{ccc}
    G\longrightarrow H& \Longrightarrow & G^{^{r\over s}}\longrightarrow H^{^{{r\over s}}}. \\
  \end{array}
$$
 \end{alphlem}
  \begin{alphlem}\label{haal}{\rm\cite{MR2587748}}
Let $r$, $s$, $p$, and $q$ be non-negative integers and G be a
graph. Then
$$(G^{^{2r+1\over 2s+1}})^{2p+1\over 2q+1}\longrightarrow G^{^{(2r+1)(2p+1)\over (2s+1)(2q+1)}}.$$
 \end{alphlem}
For a given graph $G$  with $v\in V(G)$, set
$$N_i(v)= \{u|\ {\rm there\ is\ a\ walk\ of\ length}\ i\ {\rm joining}\
u\ {\rm and}\ v\}.$$

For two subsets $A$ and $B$ of the vertex set of a graph $G$, we
write $A \bowtie B$ if every vertex of $A$ is joined to every
vertex of $B$. Also, for any non-negative integer $s$, define the
graph $G^{^{\stackrel{1}{\widetilde{2s+1}}}}$ as follows.
$$V(G^{^{\stackrel{1}{\widetilde{2s+1}}}})= \{(A_1,\ldots,A_{s+1})|\ A_i\subseteq V(G), |A_1|=1,
\varnothing\not =A_i\subseteq N_{i-1}(A_1) ,\,i\leq s+1\}.$$

Two vertices $(A_1,\ldots,A_{s+1})$ and $(B_1,\ldots,B_{s+1})$ are
adjacent in $G^{^{\stackrel{1}{\widetilde{2s+1}}}}$ if for any
$1\leq i\leq s$ and $1\leq j\leq s+1$, $A_i\subseteq B_{i+1}$,
$B_i\subseteq A_{i+1}$, and $A_{j}\bowtie B_{j}$. Here is
the definition of dual power as a functor as follows. Let $r$ and $s$ be
non-negative integers. For any graph $G$ define the graph
$G^{^{\stackrel{2r+1}{\widetilde{2s+1}}}}$ as follows
$$G^{^{\stackrel{2r+1}{\widetilde{2s+1}}}}= \left(G^{^{\stackrel{1}{\widetilde{2s+1}}}}\right)^{2r+1}.$$
These powers, in sense of graph homomorphism, inherit several
properties from power in numbers.

\begin{alphlem}\label{equi}{\rm\cite{MR2587748}}
Let $r$, $p$, and $q$ be non-negative integers. For any graph $G$
we have
\begin{description}
\item[a)] $G^{^{(2r+1)(2p+1)\over (2r+1)(2q+1)}}\longleftrightarrow G^{^{2p+1\over 2q+1}}$.
\item[b)] $G^{^{\stackrel{(2r+1)(2p+1)}{\widetilde{(2r+1)(2q+1)}}}}\longleftrightarrow G^{^{\stackrel{2p+1}{\widetilde{2q+1}}}}$.
\end{description}
\end{alphlem}

It was proved in \cite{MR2587748} that these two powers are dual
of each other as follows.
\begin{alphthm}\label{dual}{\rm\cite{MR2587748}}
Let G and H be two graphs. Also, assume that ${2r+1\over
2s+1} < {\rm og}(G)$ and $2s + 1 <{\rm og}(H^{^{\stackrel{1}{\widetilde{2r+1}}}})$.
 We have $$G^{^{2r+1\over2s+1}}\longrightarrow H \Longleftrightarrow G\longrightarrow H^{^{\stackrel{2s+1}{\widetilde{2r+1}}}}.$$
\end{alphthm}
Now, we consider the parameter $\theta_i(G)$ which in some sense measures the homomorphism capabilities of $G$.
\begin{defin}{ Assume that $G$ is a non-bipartite
graph. Also, let $i \geq -\chi(G)+3$ be an integer. The {\it $i$th
power thickness} of $G$ is defined as follows.
$$\theta_i(G) = \sup\{{2r+1 \over 2s+1}| \chi(G^{{2r+1 \over 2s+1}})\leq \chi(G)+i,
{2r+1 \over 2s+1}< {\rm og}(G) \}.$$ For simplicity, when i = 0, the
parameter is called the power thickness of $G$ and is denoted by
$\theta(G)$. Also, when $i=\chi(G)-3$, we set $\theta_{3-\chi(G)}(G)=\mu(G)$.
}
\end{defin}
\begin{alphlem}{\rm\cite{MR2587748}}\label{NOHOM}
Let $G$ and $H$  be two non-bipartite graphs with
$\chi(G)=\chi(H)-j,\ j\geq 0$. If $G \longrightarrow H$ and $i+j
\geq -\chi(G)+3$, then
$$\theta_{i+j}(G) \geq \theta_i(H).$$
\end{alphlem}
It is interesting that $\mu(G)$ is  computed in terms of circular
chromatic number. Hence, $\theta_i(G)$'s can be considered as a
generalization of circular chromatic number.
\begin{alphthm}\label{ptcir}{\rm\cite{MR2587748}} Let  $G$ be a non-bipartite graph. Then
$$\mu(G)= {\chi_c(G)\over 3(\chi_c(G)-2)}.$$
\end{alphthm}
\section{Circular Chromatic Number of Graph Powers}
Some properties of graph powers and its close relationship to
the  circular chromatic number of non-bipartite graphs have been
studied in \cite{MR2587748}. In particular, an equivalent
definition of the circular chromatic number in terms of graph
powers was introduced as follows.
\begin{alphthm}{\rm\cite{MR2587748}}{\label{haal1}}
Let G be a non-bipartite graph with chromatic number $\chi(G)$.
\begin{description}
\item[a)] If $0 < {{2r+1}\over {2s+1}} \leq {\frac{\chi(G)}{3(\chi(G)-2)}}$,
then $\chi(G^{{2r+1}\over {2s+1}})=3$. Furthermore,
$\chi(G)\neq\chi_c(G)$ if and only if there exists a rational
number ${{2r+1}\over {2s+1}}>{\frac{\chi(G)}{3(\chi(G)-2)}}$ for
which $\chi(G^{{2r+1}\over {2s+1}})= 3$.

\item[b)] $\chi_c(G)=\inf
\{{2n+1\over n-t}| \chi(G^{2n+1\over 3(2t+1)})=3, n> t >0\}.$
\end{description}
\end{alphthm}

Here, we show that if $2r+1 < {\rm og}(K_{n\over d})$, then $K_{n\over
d}^{2r+1}$ is isomorphic to a circular complete graph.

\begin{lem}\label{cirpow}
Let $n$ and $d$ be positive integers, where $n> 2d$.
\begin{description}
\item [a)] If $r$ is a non-negative integer and ${n\over d}<{2r+1\over r}$, then
$K_{n\over d}^{2r+1}\cong K_{n\over (2r+1)d-rn}$.
\item[b)] If $s$ is a nonnegative integer, then $K_{n\over d}\longleftrightarrow K_{(2s+1)n\over sn+d}^{^{2s+1}}.$
\end{description}
\end{lem}
\begin{proof}{Let $t\leq r$ be a non-negative integer. If $i$ is an arbitrary vertex of $K_{n\over d}$, it is~not hard
to check that $N_{2t+1}(i)=\{i+(2t+1)d-tn, i+(2t+1)d-tn+1,\ldots,i-(2t+1)d+tn+1\}$, where the summation is modulo $n$.
Therefore, $K_{n\over d}^{2r+1}$
is isomorphic to the circular complete graph $K_{n\over (2r+1)d-rn}$.
The next part is an immediate consequence of part (a).
}
\end{proof}
Now, we introduce an upper bound for the circular chromatic number of graph powers.
\begin{thm}
Let $r$ and $s$ be non-negative integers and $G$ be a non-bipartite graph with circular chromatic number $\chi_c(G)$. If ${2r+1\over 2s+1}<{\chi_c(G)\over \chi_c(G)-2}$, then
$$\chi_c(G^{^{2r+1\over 2s+1}})\leq{(2s+1)\chi_c(G)\over (s-r)\chi_c(G)+2r+1}.$$
\end{thm}
\begin{proof}{
Let $\chi_c(G)= {n\over d }$. It is easy to see that if ${2r+1\over 2s+1}<{{n\over d }\over {n\over d }-2}$, then ${(2s+1)n\over sn+d}<{2r+1\over r}$.

$$
  \begin{array}{ccll}
   \vspace{.3cm}G\longrightarrow K_{n\over d}& \Longrightarrow &G^{^{2r+1\over 2s+1}}\longrightarrow (K_{n\over d})^{2r+1\over 2s+1} &({\rm By~Lemma~\ref{hada}}) \\
 ~& \Longrightarrow &G^{^{2r+1\over 2s+1}} \longrightarrow (K_{(2s+1)n\over sn+d}^{^{2s+1}})^{2r+1\over 2s+1}&({\rm By~Lemma~\ref{cirpow}(b)}) \\
        \end{array}
$$
$$
  \begin{array}{ccll}
    \vspace{.3cm} ~& \Longrightarrow &G^{^{2r+1\over 2s+1}} \longrightarrow K_{(2s+1)n\over sn+d}^{^{2r+1}}&({\rm By~Lemmas~\ref{haal}~and~\ref{equi})}\\
 \vspace{.3cm} ~& \Longrightarrow &G^{^{2r+1\over 2s+1}} \longrightarrow K_{(2s+1)n\over (2r+1)(sn+d)-r(2s+1)n}&({\rm By~Lemma~\ref{cirpow}(a))}\\
    \vspace{.3cm} ~& \Longrightarrow &{\chi_c(G^{^{2r+1\over 2s+1}})}\leq {(2s+1)n\over (s-r)n+(2r+1)d} & \\
              ~& \Longrightarrow &{\chi_c(G^{^{2r+1\over 2s+1}})}\leq  {(2s+1)\chi_c(G)\over (s-r)\chi_c(G)+(2r+1)}. & \\
  \end{array}$$
}
\end{proof}
Tardif {\rm\cite{MR2171371}} has shown that the cube root, in
sense of dual power, of any circular complete graph with circular
chromatic number less than 4, is homomorphically equivalent to a
circular complete graph.

\begin{alphlem} \label{tardif}{\rm\cite{MR2171371}} Let $n$ and $d$ be positive integers, where $n > 2d$. If ${n\over d}<4, then$
$K_{n\over
d}^{^{\stackrel{1}{\widetilde{~3~}}}}\longleftrightarrow
K_{3n\over n+d}$.
\end{alphlem}

Here is a generalization of Lemma \ref{tardif}.

\begin{lem}\label{tar1} Let $n$ and $d$ be positive integers, where $n > 2d$.
If ${n\over d}<4, then$
$$K_{n\over d}^{^{\stackrel{1}{\widetilde{2r+1}}}}\longleftrightarrow K_{(2r+1)n\over rn+d}.$$
\end{lem}
\begin{proof}{
Theorem~\ref{dual} implies that $K_{(2r+1)n\over
rn+d}\longrightarrow K_{n\over
d}^{^{\stackrel{1}{\widetilde{2r+1}}}}$ if and only if
$K_{(2r+1)n\over rn+d}^{^{2r+1}}\longrightarrow K_{n\over d}$. On
the other hand, Lemma~\ref{cirpow}(b) shows that the circular
complete graphs $K_{(2r+1)n\over rn+d}^{^{2r+1}}$ and $K_{n\over
d}$ are homomorphically equivalent. Conversely,  it is sufficient
to prove that
$$\chi_c((K_{n\over d})^{^{\stackrel{1}{\widetilde{2r+1}}}})\leq {(2r+1)n\over rn+d}.$$
Take a rational number ${2k+1\over 3^i}$ such that $ {1\over
2r+1}\leq {2k+1\over 3^i}<{1\over 2}$. It is easy to see that
$$K_{n\over d}^{^{\stackrel{1}{\widetilde{2r+1}}}}\longrightarrow (K_{n\over d})^{\stackrel{2k+1}{\widetilde{~{~3^i~}~}}}.$$
If $G$  is non-bipartite graph,  Theorem~\ref{dual} and
Lemma~\ref{equi} yield that
$G^{^{\stackrel{1}{\widetilde{~{~3^i~}~}}}}\longrightarrow
(G^{^{\stackrel{1}{\widetilde{~{~3~}~}}}})^{^{\stackrel{1}{\widetilde{~{~3^{i-1}~}~}}}}$.
Since ${n\over d}<4$, by induction on $i$ and Lemma~\ref{tardif}
we have $K_{n\over
d}^{\stackrel{2k+1}{\widetilde{~{~3^i~}~}}}\longleftrightarrow
K_{3^in\over{{3^i-1\over 2}n+d}}^{^{2k+1}} $. Therefore, there is
a homomorphism from $K_{n\over
d}^{^{\stackrel{1}{\widetilde{2r+1}}}}$ to
$K_{3^in\over{{3^i-1\over 2}n+d}}^{^{2k+1}} $. By
Lemma~\ref{cirpow}(a), two graphs $K_{3^in\over{{3^i-1\over
2}n+d}}^{^{2k+1}}$ and $K_{3^in\over{(2k+1)({{3^i-1\over
2}n+d})-k3^in}}$ are homomorphically equivalent. Hence,
$$\chi_c(K_{n\over d}^{^{\stackrel{1}{\widetilde{2r+1}}}})\leq {3^in\over{(2k+1)d+({3^i-1\over 2}-k)n}}.$$
Since the set of parameters $\{{2k+1\over 3^i}\mid k\geq 1, i\geq 1 \}$ is dense in the interval $(0,+\infty)$,
$$\chi_c(K_{n\over d}^{^{\stackrel{1}{\widetilde{2r+1}}}})\leq \inf\left\{{3^in\over{(2k+1)d+({3^i-1\over 2}-k)n}}
{\bigg |}{{1\over 2r+1}\leq {2k+1\over 3^i}}<{1\over 2}\right\}.$$
This infimum is equal to  ${(2r+1)n\over rn+d}$, as desired. }
\end{proof}

Here, we determine the circular chromatic number of some graph
powers.

\begin{thm}\label{cirsub}
Let $G$ be a non-bipartite graph with circular chromatic number
$\chi_c(G)$. Also, assume that $r$ and $s$ are non-negative
integers. Then we have $\chi_c(G^{^{1\over
2s+1}})=\frac{(2s+1)\chi_c(G)}{s\chi_c(G)+1}.$ Moreover, If
$\chi_c(G^{^{2r+1\over 2s+1}})<4$, then
$$\chi_c(G^{^{2r+1\over 2s+1}})={(2s+1)\chi_c(G)\over (s-r)\chi_c(G)+2r+1}.$$
\end{thm}
\begin{proof}{
Note that, in view of Theorem~\ref{dual},  $G^{^{1\over
2s+1}}\longrightarrow K_{(2s+1)\chi(G) \over s\chi(G)+1}$ if and
only if $G\longrightarrow K_{(2s+1)\chi(G) \over
s\chi(G)+1}^{2s+1}$. On the other hand, by using
Lemma~\ref{cirpow}(a), two graphs $K_{(2s+1)\chi(G) \over
s\chi(G)+1}^{2s+1}$ and  $K_{(2s+1)\chi(G) \over (2s+1)}$ are
homomorphically equivalent. Consequently, $\chi_c(G^{^{1\over
2s+1}})<{2s+1\over s}$. Let ${n\over d}<{2s+1\over s}$.
$$
  \begin{array}{ccll}
  \vspace{.3cm}\chi_c(G^{^{1\over 2s+1}})\leq {n\over d}   & \Longleftrightarrow & G^{^{1\over 2s+1}}\longrightarrow K_{n\over d}& \\
  \vspace{.3cm} ~& \Longleftrightarrow & G\longrightarrow K_{n\over d}^{^{2s+1}}&({\rm By~Theorem~\ref{dual}}) \\
     \vspace{.3cm} ~& \Longleftrightarrow & G\longrightarrow K_{n\over (2s+1)d-sn}&({\rm By~Lemma~\ref{cirpow}}({\rm a})) \\
    \vspace{.3cm} ~& \Longleftrightarrow & \chi_c(G)\leq \chi_c(K_{n\over (2s+1)d-sn})& \\
       ~& \Longleftrightarrow &{(2s+1)\chi_c(G)\over s\chi_c(G)+1}\leq {n\over d} & \\
  \end{array}
$$

To prove the next part, it suffices to show that for any
$2\leq{n\over d}<4$, $\chi_c(G^{^{2r+1\over 2s+1}})\leq {n\over
d}$ is equivalent to ${(2s+1)\chi_c(G)\over
(s-r)\chi_c(G)+2r+1}\leq {n\over d}$.
 Assume that $\chi_c(G^{^{2r+1\over 2s+1}})\leq {n\over d}<4$.

$$
  \begin{array}{ccll}
  \vspace{.3cm}\chi_c(G^{^{2r+1\over 2s+1}})\leq {n\over d}   & \Longleftrightarrow &
  G^{^{2r+1\over 2s+1}}\longrightarrow K_{n\over d}& \\
   \vspace{.3cm} ~& \Longleftrightarrow & G^{^{1\over 2s+1}}\longrightarrow
   K_{n\over d}^{^{\stackrel{1}{\widetilde{2r+1}}}}&({\rm By~Theorem~\ref{dual}}) \\
   \vspace{.3cm} ~& \Longleftrightarrow & G^{^{1\over 2s+1}}\longrightarrow
   K_{(2r+1)n\over rn+d}&({\rm By~Lemma~\ref{tar1}}) \\
    \vspace{.3cm} ~& \Longleftrightarrow & \chi_c(G^{^{1\over 2s+1}})\leq {(2r+1)n\over rn+d}& \\
    \vspace{.3cm} ~& \Longleftrightarrow & {(2s+1) \chi_c(G)\over s\chi_c(G)+1}\leq{{(2r+1){n\over d}}\over
    r{n\over d}+1}& \\
       ~& \Longleftrightarrow &{(2s+1)\chi_c(G)\over (s-r)\chi_c(G)+2r+1}\leq {n\over d}. & \\
  \end{array}
$$
}
\end{proof}

\begin{cor}
Let $r$ and $s$ be non-negative integers and $G$ be a
non-bipartite graph. if ${2r+1 \over 2s+1} \leq {\chi_c(G) \over 3(\chi_c(G)-2)}$,
then $\chi_c(G^{^{\frac{2r+1}{2s+1}}})=\frac{(2s+1)\chi_c(G)}{(s-r)\chi_c(G)+2r+1}$.
\end{cor}
\begin{proof}{
Since ${2r+1 \over 2s+1} \leq {\chi_c(G) \over 3(\chi_c(G)-2)}$, Theorem \ref{haal1} implies that
$\chi_c(G^{^{\frac{2r+1}{2s+1}}})\leq 3$. Now, by the previous theorem, we have $\chi_c(G^{^{\frac{2r+1}{2s+1}}})=\frac{(2s+1)\chi_c(G)}{(s-r)\chi_c(G)+2r+1}$.
}
\end{proof}
\begin{cor}
Let $r$ and $s$ be non-negative integers and $G$ be a
non-bipartite graph such that $\chi_c(G^{^{2r+1\over 2s+1}})<4$.
Then we have
$$\mu(G^{^{2r+1\over 2s+1}})={2s+1\over 2r+1}\mu(G)={2s+1\over 3(2r+1)}{\chi_c(G)\over (\chi_c(G)-2)}.$$
\end{cor}
Given a rational number ${n \over d}$, a rational number ${n'
\over d'}$ is unavoidable by ${n \over d}$ if every graph $G$
with $\chi_c(G) = {n \over d}$ contains a subgraph $H$ with
$\chi_c(H) = {n' \over d'}$. It is known \cite{MR2003514} if $m$
is an integer and $m < {n \over d}$, then $m$ is unavoidable by
${n\over d}$.

Suppose $(n, d)=1$, i.e., $n$ and $d$ are coprime. Let $n'$ and
$d'$ be the unique integers such that $0 < n' < n$ and
$nd'-n'd=1$. We call ${n' \over d'}$ the lower parent of ${n
\over d}$, and denote it by $F({n \over d})$. The following
question was posed in {\rm \cite{MR1815614,MR2249284}}.

\begin{alphqu}\label{zhuqu}{\rm \cite{MR1815614,MR2249284}}
Is true that for every rational ${n \over d}>2$, $F({n \over d})$
is unavoidable by ${n \over d}${\rm ?}
\end{alphqu}

Here, we give a negative answer to the aforementioned question.

\begin{cor}
Let $k$ be a positive integer. Then there exists a graph $G$ with
$\chi_c(G)={9k+3\over 3k+2}$ such that $G$ does not contain any
subgraph with circular chromatic number  equal to ${6k+1\over
2k+1}$.
\end{cor}
\begin{proof}{
Let $n=9k+3$, $d=3k+2$, $n'=6k+1$, and $d'=2k+1$. Obviously,
$nd'-n'd=1$. By Theorem~\ref{cirsub}, we have
$\chi_c(K_{3k+1}^{^{1\over3}})={9k+3\over 3k+2}$. Suppose that
$e\in E(K_{3k+1}^{^{1\over3}})$. It is readily seen that there
exists a homomorphism from $K_{3k+1}^{^{1\over3}}\setminus e$ to
$K_{3k}^{^{1\over3}}$. Hence, if $H$ is a proper subgraph of
$K_{3k+1}^{^{1\over3}}$, then  $\chi_c(H)\leq
\chi_c(K_{3k}^{^{1\over3}})={9k\over 3k+1}< {6k+1\over 2k+1}$.
Therefore, $G$ contains no subgraph  with circular chromatic
number $n'\over d'$.
}
\end{proof}

It should be noted that one can introduce more rational numbers
 such that their lower parents are not  unavoidable. For instance,
 we show that ${15n+7 \over 5n+4}$ is not unavoidable by ${18n+9 \over
 6n+5}$. To see this, for $d\geq 2$ and $n\geq 3$, define the graph $H_d(K_n)$ as follows. Let
$G_1,\ldots ,G_d$ be $d$ graphs such that each of them is
isomorphic to the complete graph $K_n$. Assume that $v_iw_i \in
E(G_i)$ for any $1\leq i \leq d$.  The graph $H_d(K_n)$ obtained
from the disjoint union of $G_1\cup \cdots \cup G_d$ by
identifying the vertices $w_i$ with $v_{i+1}$ for any $1\leq i
\leq d-1$, deleting the edges $v_iw_i$ for any $1\leq i \leq d$,
and by adding the edge $v_1w_d$. In fact, it is a simple matter to
check that $H_d(K_n)$ follows by applying Haj\'{o}s construction
to the complete graphs $G_1,\ldots ,G_d$. Hence,
$\chi(H_d(K_n))=n$ and the graph $H_d(K_n)$ is a critical graph,
i.e., $\chi(H_d(K_n)\setminus e)=n-1$ for any $e\in E(H_d(K_n))$.

Now, we show that $\chi_c(H_d(K_n))={d(n-1)+1 \over d}$. To see
this, assume that $V(G_i)=\{v_i,u_{i2}, \ldots, u_{i(n-1)},
w_i\}$. Define a coloring $c:V(H_d(K_n)) \rightarrow
\{1,2,\ldots,dn-d+1\}$ as follows. For any $1\leq i\leq d$ and
$2\leq j\leq n-1$, set $c(u_{ij})=(j-1)d+i$,  $c(w_i)=i$, and
$c(v_1)=d(n-1)+1$. It is easy to check that $c$ is a
$(d(n-1)+1,d)$ coloring of $H_d(K_n)$.

On the other hand, it is straightforward to check that the
independence number of $H_d(K_n)$ is equal to $d$.
Consequently, $\chi_c(H_d(K_n))={dn-d+1 \over d}$.

The graph $H_2(K_{3n+2})^{^{1\over3}}$ has circular chromatic number $\frac{18n+9}{6n+5}$.
 It is readily seen that there is a homomorphism from
$H_2(K_{3n+2})^{^{1\over3}}\setminus e$ to $K_{3n+1}^{^{1\over3}}$.
Hence, if $H$ is a proper subgraph of $H_2(K_{3n+2})^{^{1\over3}}$,
then $\chi_c(H)\leq \chi_c(K_{3n+1}^{^{1\over3}})={9n+3\over 3n+2}<
{15n+7\over 5n+4}$. Therefore, $H_2(K_{3n+2})^{^{1\over3}}$ contains no subgraph
with circular chromatic number ${15n+7\over 5n+4}$.

Let $\zeta(G)$ be the minimum number of vertices of $G$, necessary
to be deleted, in order to reduce the chromatic number of the
graph.

\begin{alphqu}{\rm\cite{MR2340388}}{\label{zeta}} Let $\chi_{_{c}}(G)=\frac{n}{d}$, where $(n,d)=1$ and
$n=(\chi(G)-1)d+r$. Is it true that $\zeta(G) \geq r$ {\rm ?}
\end{alphqu}

When $G$ is a critical graph, we have $\zeta(G)=1$. If the
aforementioned question is true, then for every critical graph $G$
with $\chi(G)=n$, its circular chromatic number is equal to
${dn-d+1 \over d}$ for an appropriate $d$. It is worth noting
that $H_d(K_n)$ is a critical graph with
$\chi_c(H_d(K_n))={dn-d+1 \over d}$.
\section{Fractional and Multichromatic Number}
As usual, we denote by $[m]$ the set $\{1, 2, \ldots, m\}$, and
denote by ${[m] \choose n}$ the collection of all $n$-subsets of
$[m]$. The {\em Kneser graph} ${\rm KG}(m,n)$ ({\em resp. the
generalized Kneser graph} ${\rm KG}(m,n,s)$) is the graph on the
vertex set ${[m] \choose n}$, in which two distinct vertices $A$
and $B$ are adjacent if and only if $A \cap B = \varnothing$
(resp. $|A\cap B|\leq s$). It was conjectured by Kneser
\cite{Kneser} in 1955, and proved by Lov\'{a}sz \cite{MR514625}
in 1978, that $\chi({\rm KG}(m,n))=m-2n+2$. The fractional
chromatic number is defined as a generalization of the chromatic
number as follows
$$\chi_f(G)=\inf\{{m\over n}| G\rightarrow {\rm KG}(m,n)\}.$$

An $n$-tuple coloring of graph $G$ with $m$ colors assigns to each
vertex of $G$, an $n$-subset of $[m]$ so that adjacent vertices
receive disjoint sets. Equivalently, $G$ has an $n$-tuple
coloring with $m$ colors if there exists a homomorphism from $G$
to ${\rm KG}(m,n)$. The $n$th multichromatic number of $G$,
denoted by $\chi_n(G)$, is the smallest $m$ such that $G$ has a
$n$-tuple coloring with $m$ colors. These colorings were first
studied in the early 1970s and the readers are referred to
\cite{MR1475894,MR1481157,MR1614286} for more information.

\begin{alphthm}{\rm\cite{MR1883597}}\label{pir}
Suppose that $m$ and $n$ are positive integers with $m>2n$. Then
the following two conditions on non-negative integers $k$ and $l$
are equivalent.
\begin{itemize}
\item For any two {\rm(}not necessarily distinct{\rm)} vertices $A$ and $B$ of ${\rm KG}(m,n)$ with
$|A\cap B|=k$, there is a walk of length exactly $l$ in ${\rm KG}(m,n)$ beginning at $A$ and ending at $B$.
\item $l$ is even and $k\geq n-{l\over 2}(m-2n)$, or $l$ is odd and $k\leq {l-1\over 2}(m-2n)$.
\end{itemize}\end{alphthm}

In view of Theorem~\ref{cirsub}, we have $\chi_c(G^{^{1\over
2s+1}})={{(2s+1)\chi_c(G)}\over{s\chi_c(G)+1}}$. Here, we present
a tight upper bound for the fractional chromatic number of subdivision
graphs.
\begin{thm}\label{fracb}
Let $G$ be a non-bipartite graph and $s$ be a non-negative
integer. Then
$$\chi_f(G^{^{1\over
2s+1}})\leq{{(2s+1)\chi_f(G)}\over{s\chi_f(G)+1}}.$$
\end{thm}
\begin{proof}{
Let $f$ be a homomorphism from $G$ to ${\rm KG}(m,n)$.
We claim that there is a homomorphism from $G$ to the generalized
Kneser graph ${\rm KG}((2s+1)m,sm+n,(m-2n)s)$. To see this, for
every vertex $v\in V(G)$, define $g(v)$ as follows
$$
\displaystyle \bigcup_{i \not\in f(v)}\{(i-1)(2s+1)+1,\ldots,(i-1)(2s+1)+s\}
\bigcup_{i \in f(v)}\{(i-1)(2s+1)+s+1,\ldots,i(2s+1)\}.
$$
It is easy to see that, for any vertex $v\in V(G)$, $|g(v)|=sm+n$.
Also, if $u$ and $v$ are two adjacent vertices in $G$, then
$|g(u)\cap g(v)|=(m-2n)s$. Now, in view of Theorem~\ref{pir}, we
have
 $${\rm KG}((2s+1)m,sm+n,(m-2n)s)\longleftrightarrow {\rm KG}((2s+1)m,sm+n)^{^{2s+1}}.$$
Let $G\longrightarrow {\rm KG}(m,n)$ and $\chi_f(G)={m\over n}$.
By the previous discussion, there is a homomorphism from $G$ to
${\rm KG}((2s+1)m,sm+n,(m-2n)s)$. Now, Theorem~\ref{dual} implies that
 $G^{^{1\over2s+1}} \longrightarrow{\rm KG}((2s+1)m,sm+n).$ Hence, $\chi_f(G^{^{1\over2s+1}} )\leq
{{(2s+1)\chi_f(G)}\over{s\chi_f(G)+1}}.$ }
\end{proof}

Equality does~not always hold in Theorem~\ref{fracb}. For
instance, consider the graph $K_{10}^{1\over 3}$. We know that the
third power of the Petersen graph $P^3$ is isomorphic to
$K_{10}$. Hence, in view of Lemma~\ref{haal}, there exists a
homomorphism from $K_{10}^{1\over 3}$ to the Petersen graph.
Consequently, $\chi_f(K_{10}^{1\over 3})\leq {5 \over 2}$ which is
less than ${31 \over 11}$.

It is simple to see that there exists a homomorphism from
$G^{^{1\over 2n+1}}$ to $C_{2n+1}$. On the other hand, the odd
cycle $C_{2n+1}$ is an induced subgraph of the Kneser graph
$KG(2n+1,n)$. Therefore,  if $G$ is a non-bipartite graph and
$s\geq n$, then $\chi_n(G^{^{1\over 2s+1}})=2n+1$.
\begin{thm}\label{multi}
Let $G$ be a non-bipartite graph. If  $i,n$ and $s$ are positive integers such that $is=n-1$, then
$$\chi_n(G^{^{1\over 2s+1}})\leq 2n+i\Longleftrightarrow \chi(G)\leq {2n+i\choose n}$$
\end{thm}
\begin{proof}{
$$
  \begin{array}{ccll}
 \vspace{.3cm} \chi_n(G^{^{1\over 2s+1}})\leq 2n+i  & \Longleftrightarrow &
  G^{^{1\over 2s+1}}\longrightarrow {\rm KG}(2n+i,n) & \\
   \vspace{.3cm} ~& \Longleftrightarrow & G\longrightarrow
   {\rm KG}(2n+i,n)^{2s+1}&({\rm By~Theorem~\ref{dual}}) \\
   \vspace{.3cm} ~& \Longleftrightarrow & G\longrightarrow
   {\rm KG}(2n+i,n,is)&({\rm By~Theorem~\ref{pir}}) \\
   \vspace{.3cm} ~& \Longleftrightarrow & G\longrightarrow
   {\rm KG}(2n+i,n,n-1)& \\
\vspace{.3cm} ~& \Longleftrightarrow & G\longrightarrow
   K_{2n+i\choose n}& \\
    \vspace{.3cm} ~& \Longleftrightarrow & \chi(G)\leq {2n+i\choose n}. \\
  \end{array}
$$
}
\end{proof}
We know that $\chi_2(G^{^{1\over 2s+1}})=5$ whenever $s\geq 2$.
The following corollary, which is an immediate consequence of the
aforementioned theorem, determines the other cases.

\begin{cor}
Let $G$ be a non-bipartite graph. If $\chi(G)\leq 10$, then
$\chi_2(G^{^{1\over 3}})=5$. Otherwise, $\chi_2(G^{^{1\over
3}})=6$.
\end{cor}

\end{document}